\documentclass[11pt]{amsart}
\usepackage{amsmath}
\usepackage{amsfonts}
\usepackage{amssymb}

\theoremstyle{plain}
\newtheorem{thm}{Theorem}[section]
\newtheorem{cor}[thm]{Corollary}

\newtheorem{lemma}[thm]{Lemma}
\newtheorem{prop}[thm]{Proposition}
\newtheorem{repthm}{Theorem}

\theoremstyle{definition}

\newtheorem{prob}[thm]{Problem}

\newcommand{\bb}[1]{\mathbb{#1}}
\newcommand{\cl}[1]{\mathcal{#1}}
\newcommand{\ff}[1]{\mathfrak{#1}}
\newcommand{\mr}[1]{\mathrm{#1}}

\newcommand{\norm}[1]{\left\| #1 \right\|}
\newcommand{\inp}[2]{\left\langle #1 , #2 \right\rangle}

\newcommand{\ip}[1]{\left\langle #1 \right\rangle}
\newcommand{\dist}{\operatorname{dist}}
\newcommand{\ol}{\overline}

\renewcommand{\phi}{\varphi}
\newcommand{\ca}{\mathrm{C}^*}
\newcommand{\ce}{\ca_e}
\newcommand{\diag}{\operatorname{diag}}
\newcommand{\spn}{\operatorname{span}}

\newcommand{\qand}{\quad\text{and}\quad}

\newcommand{\qfor}{\quad\text{for}\quad}
\newcommand{\qforal}{\quad\text{for all}\quad}
\newcommand{\AND}{\text{ and }}

\begin{document}

\title[A Constrained Nevanlinna-Pick Interpolation Problem]{A Constrained Nevanlinna-Pick\\ Interpolation Problem}

\author[K.R.~Davidson]{Kenneth R. Davidson}
\thanks{First author supported in part by a grant from NSERC}
\address{Pure Mathematics Department, University of Waterloo, 
Waterloo, Ontario N2L 3G1, Canada}
\email{krdavids@math.uwaterloo.ca}

\author[V.I.~Paulsen]{Vern I.~Paulsen}
\address{Department of Mathematics, University of Houston,
Houston, Texas 77204-3476, U.S.A.}
\email{vern@math.uh.edu}
\thanks{Second and third authors supported in part by a grant from NSF}

\author[M.~Raghupathi]{Mrinal Raghupathi}
\address{Department of Mathematics, University of Houston,
Houston, Texas 77204-3476, U.S.A.}
\email{mrinal@math.uh.edu}

\author[D.~Singh]{Dinesh Singh}
\address{Department of Mathematics, University of Delhi, New Delhi,
  India}
\email{dineshsingh1@gmail.com}
\thanks{Fourth author was a Visiting Professor at U. Houston during part of this period}

\begin{abstract}
We obtain necessary and sufficient conditions for Nevanlinna-Pick
 interpolation on the unit disk with the additional restriction
that all analytic interpolating functions satisfy $f'(0)=0.$
Alternatively, these results can be interpreted as interpolation results 
for $H^{\infty}(V),$ where $V$ is the intersection of the bidisk 
with an algebraic variety. 
We use an analysis of C*-envelopes to show that these same conditions
do not suffice for matrix interpolation.
\end{abstract}

\maketitle

\section{Introduction} \label{S:intro}

The classic Nevanlinna-Pick interpolation result says that given
$n$ distinct points, $z_1,\dots,z_n,$ in the open unit disk, $\bb
D,$ and $n$ complex numbers, $w_1,\dots,w_n$ and $A>0,$ then there exists
an analytic function $f$ on $\bb D$ with $f(z_i)=w_i$ for $i=1,\dots,n$ and
$\|f\|_{\infty} \le A$ if and only if the $n \times n$ matrix, 
\[
 \left[ \frac{A^2- w_i \overline{w_j}}{1-z_i \overline{z_j}} \right]
\]
is positive semidefinite, where $\|f\|_{\infty} := \sup \{ |f(z)|: z \in \bb D \}.$

In this paper we give two very different sets of necessary and
sufficient conditions for the classical Nevanlinna-Pick problem with
one additional constraint. Namely, that $f^{\prime}(0)=0.$  
We define the algebra
\[
 H^{\infty}_1 = H^{\infty}_1(\bb D) :=
  \{ f \in H^{\infty}(\bb D): f^{\prime}(0)=0 \},
\]
so that our constraint is simply the requirement that functions belong to this algebra.

Our main result is analogous to M.B.~Abrahamse's \cite{Ab} interpolation results for finitely
connected domains. If $\cl R$ is a bounded domain in the complex plane
whose boundary consisted of $p+1$ disjoint analytic Jordan curves,
then Abrahamse identified
a family of reproducing kernel Hilbert spaces $H^2_{\alpha}(\cl R)$
indexed by $\alpha$ in the $p$-torus $\bb T^p$,
with corresponding kernels $K_{\alpha}(z,w)$. He proved that if
$z_1,\dots,z_n$ are $n$ distinct points in $\cl R$ and $w_1,\dots,w_n$ are
complex numbers, then there exists an analytic function $f$ on $\cl R$
such that $f(z_i)=w_i$ and $\sup \{ |f(z)| : z \in \cl R \} \le A$ if
and only if 
\[
 \Big[ (A^2 -w_i\overline{w_j})K_{\alpha}(z_i,z_j) \Big] \ge 0
 \qforal \alpha \in \bb T^p.
\]

In a similar fashion, we identify a family of reproducing kernel
Hilbert spaces of analytic functions on $\bb D$, denoted $H^2_{\alpha,
  \beta}(\bb D)$, indexed by points on the sphere in complex 2-space,
$|\alpha|^2+|\beta|^2=1$. The corresponding kernels are given by
\[
 K^{\alpha,\beta}(z,w) 
 = (\alpha + \beta z)\overline{(\alpha + \beta w)} 
 + \frac{z^2 \overline{w^2}}{1-z \overline{w}}.
\]
We will prove that these kernel functions play a similar role for our constrained 
interpolation problem to the role played by Abrahamse's kernel functions for interpolation on finitely connected domains.

In Section \ref{S:npint}, we will prove the following theorem.

\begin{thm}\label{thm1} Let $z_1,\dots,z_n$ be distinct points in $\bb D$, 
and let $w_1,\dots,w_n$ be complex numbers. Then there exists an
analytic function $f$ on $\bb D$ with $\|f\|_{\infty} \le A$ and 
$f^{\prime}(0)=0$ such that $f(z_i)=w_i$ for $i=1,\dots,n$ if and only if 
\[
  \Big[ (A^2-w_i \overline{w_j})K^{\alpha,\beta}(z_i,z_j) \Big]
\]
is positive semidefinite for all $|\alpha|^2+|\beta|^2 =1.$
\end{thm}

Note that the kernels $K^{\alpha,\beta}$ and 
$K^{e^{i\theta}\alpha,e^{i\theta}\beta}$ coincide.
That is, this family of kernels is parameterized by
the set of complex lines in $\bb C^2.$ 
The complex projective 2-sphere will be denoted by $PS^2$.
So in Theorem \ref{thm1}, it is sufficient to consider the kernels
for $(\alpha, \beta)= (re^{i \theta}, \sqrt{1-r^2})$ for $0 \le r \le 1$ 
and $e^{i \theta} \in \bb T$;
and when $r=1$, these kernels are all equal to the kernel for $(1,0)$. 
This also shows that $PS^2$ may be identified topologically with the real 2-sphere. 
Consequently, the conditions in Theorem \ref{thm1} can be regarded as involving
parameters indexed by a real 2-sphere.

This parameterization will be used in section~\ref{S:Cenv},
where we want exactly one representative for each kernel in order to discuss representations of certain quotient algebras.

One can also prove an interpolation result
for this situation that involves an existential quantifier. Given
$\lambda \in \bb D$, let $\phi_{\lambda}(z) = \frac{z - \lambda}{1- \ol{\lambda}z}$
denote the elementary M\"{o}bius transformation
that sends $\lambda$ to $0$.

\begin{thm}\label{thm3} 
Let $z_1,\dots,z_n$ be distinct points in $\bb D$, 
and let $w_1,\dots,w_n$ be in $\bb D$. Then there exists an
analytic function $f \in H^\infty_1$ with $\|f\|_{\infty} \le 1$ 
such that $f(z_i)=w_i$ for $i=1,\dots,n$ if and only if
there exists $\lambda \in \bb D$ so that
\[ \left[ \frac{z_i^2 \overline{z_j^2} - \phi_{\lambda}(w_i)
  \overline{\phi_{\lambda}(w_j)}}{1-z_i \overline{z_j}} \right]\]
is positive semidefinite.
\end{thm}

These theorems are also in a certain sense complementary. 
If one is ``lucky'' enough to have found a point $\lambda$ 
such that the matrix of Theorem \ref{thm3} is positive semidefinite, 
then one knows that an interpolating function exists. 
On the other hand, if one is ``lucky'' enough to find a pair 
$(\alpha,\beta)$ such that the corresponding matrix in 
Theorem \ref{thm1} is not positive semidefinite, 
then one knows that the interpolation problem has no solution.
  
Abrahamse proves that the spaces $H^2_{\alpha}(\cl R)$ in some sense
parameterize all ``nice'' models for Hilbert spaces that are modules
over the algebra $H^{\infty}(\cl R)$.  His proof of the
interpolation theorem over these finitely connected domains
demonstrates why finding these models plays a central role.
In a similar fashion, we prove that the spaces 
$H^2_{\alpha , \beta}(\bb D)$ serve as models for the spaces 
which are modules over $H^{\infty}_1$.
In particular, in section 2, we prove the following analogue of the Beurling and
Helson--Lowdenslager  theorems:

\begin{thm}\label{thm2} Let $\cl M$ be a norm
  closed subspace of  $L^2(\bb T)$ which is invariant for $H^{\infty}_1(\bb D)$,
  but is not invariant for $H^{\infty}(\bb D)$.
  Then there exist $\alpha, \beta$ in $\bb C$ with $|\alpha|^2+ |\beta|^2 =1$
  and $\alpha \ne 0$ and a unimodular function, $J$, such that
  $\cl M = J H^2_{\alpha , \beta}(\bb D).$
  When $\cl M$ is a subspace of $H^2(\bb D)$, $J$ is an inner function.
\end{thm}

Theorem \ref{thm2} is a refinement of the analogue of Beurling's theorem that 
was obtained for this algebra in an earlier paper \cite{PS}.

There are several reasons for studying the function theory of the algebra $H^{\infty}_1.$ 
Agler and McCarthy\cite{AMc2} introduced the study of the algebras $H^{\infty}(V)$ of 
bounded, analytic functions on {\em embedded disks.}  
Consider the map $\phi: \bb D \to \bb D^2$ defined by $\phi(z) = (z^2,z^3)$,
and set $V= \phi(\bb D)$.
Then $V$ is the intersection of an algebraic variety with the bidisk. 
This is a very simple example of the embedded disks studied by Agler--McCarthy. 
Moreover, the isometric homomorphism induced by composition
$\phi^*: H^{\infty}(V) \to H^{\infty}(\bb D)$ can be easily seen to be an 
isometric isomorphism onto $H^{\infty}_1(\bb D),$ 
since both $H^{\infty}_1(\bb D)$ and $\phi^*(H^{\infty}(V))$ can be seen to be the unital subalgebra of $H^{\infty}(\bb D)$ generated by $z^2$ and $z^3$. 
This fact is also pointed out in \cite{PS}, and was one of the original motivations 
for obtaining a Beurling-type theorem for the algebra $H^{\infty}_1.$
The second reason is the work of Solazzo \cite{So}, who studied Nevanlinna-Pick 
interpolation with the additional constraint that all interpolating functions have 
equal values for some finite set of points, or, equivalently, for finite codimension 
subalgebras of $H^{\infty}(\bb D)$ consisting of the bounded, analytic functions 
that are constant on a given finite set. 
Thus, it was natural to extend this work by considering the simplest example of 
a subalgebra determined by a derivative condition.

In section 4, we use the interpolation results to obtain a distance formula.
In section 5, we use C*-algebra techniques to prove the failure of the matrix-valued analogues of our interpolation results. In sections 6 and 7, we discuss two-point interpolation and attempt to compute the naturally induced pseudo-hyperbolic metric arising from our algebra. Section 8 contains a summary of open questions about interpolation for this algebra.

\section{An Invariant Subspace Characterization} \label{S:invariant}

In \cite{PS}, the subspaces of the Hardy space $H^2$ 
that are invariant under multiplication by the functions 
in the algebra $H^{\infty}_1$ are characterized. 
We begin by giving a somewhat stronger version of this result with a new proof.

Given complex numbers $\alpha$ and $\beta$ with $|\alpha|^2 + |\beta|^2 =1,$ 
we let $H^2_{\alpha, \beta}$ denote the codimension one subspace of $H^2,$ 
\[ H^2_{\alpha, \beta} := \spn \{ \alpha + \beta z,\ z^2H^2 \} .\]

It is easily checked that $H^2_{\alpha, \beta}$ is invariant for $H^{\infty}_1$.
Also $H^2_{\alpha, \beta} = H^2_{\delta, \gamma}$ if and only if the 
vectors $(\alpha, \beta)$ and $(\delta, \gamma)$ are scalar multiples of each other.

When $\alpha \ne 0$, it is easy to see that $H^2_{\alpha, \beta}$ 
is {\em not} invariant under multiplication by $H^{\infty}$.
However when $\alpha = 0$,  $H^2_{\alpha, \beta} = H^2_{0,1} = z H^2$,
which is invariant under multiplication by $H^{\infty}$.

Since these are subspaces of $H^2$, they are reproducing kernel
Hilbert spaces of analytic functions on $\bb D$. 
An orthonormal basis for $H^2_{\alpha, \beta}$ is given by 
$\{ \alpha + \beta z \} \cup \{z^n: n \ge 2 \}$.
Hence the reproducing kernel for this space is given by
\begin{align*}
 K_{\alpha, \beta}(z,w) &= 
 (\alpha + \beta z) \overline{(\alpha + \beta w)} + \sum_{n \ge 2} z^n \overline{w^n} \\
 &= (\alpha + \beta z) \overline{(\alpha + \beta w)} + \frac{z^2 \overline{w^2}}{1 - z \overline{w}}.
\end{align*}
When $\alpha = 0$, this simplifies to $K_{0,1}(z,w)= \dfrac{z \overline{w}}{1- z \overline{w}}.$

We are now ready to prove:

\renewcommand{\therepthm}{\ref{thm2}}
\begin{repthm}  Let $\cl M$ be a norm
  closed subspace of  $L^2(\bb T)$ which is invariant for $H^{\infty}_1$,
  but is not invariant for $H^{\infty}$.
  Then there exist scalars $\alpha, \beta$ in $\bb C$ with $|\alpha|^2+ |\beta|^2 =1$
  and $\alpha \ne 0$ and a unimodular function, $J$, such that
  $\cl M = J H^2_{\alpha , \beta}$. 
  When $\cl M$ is a subspace of $H^2$, $J$ is an inner function.
\end{repthm}

\begin{proof}
Let $\tilde{\cl M} = \overline{H^\infty \cdot \cl M}$.
Observe that 
\[
 \tilde{\cl M} \supset \cl M = H^\infty_1 \cdot \cl M 
 \supset \overline{z^2 H^\infty \cdot \cl M} = z^2  \tilde{\cl M} .
\]
By the Helson--Lowdenslager Theorem, the invariant subspaces of 
$H^\infty$ have the form $L^2(E)$ for some measurable subset $E \subset \bb T$
or $J H^2$ for some unimodular function $J$.

If $\tilde{\cl M} = L^2(E)$,
then $\tilde{\cl M} = z^2 \tilde{\cl M} = \cl M$ is invariant for $H^\infty$.
Thus $\tilde{\cl M} = J H^2$ for a unimodular function $J$.
So $JH^2 \supset \cl M \supset z^2J H^2$.
Both containments must be strict, as $\cl M$ is not invariant for $H^\infty$.
As $z^2 \tilde{\cl M}$ is codimension 2 in $\tilde{\cl M}$,
there is an $(\alpha, \beta)$ in the unit sphere of $\bb C^2$ so that 
\[
 \cl M = \langle J(\alpha + \beta z) \rangle \oplus z^2J H^2 = J H^2_{\alpha,\beta} .
\]
The last claim is immediate.
\end{proof}

\section{The Interpolation Theorem} \label{S:npint}

In this section we present the proofs of Theorems \ref{thm1} and \ref{thm3}.

Recall that if $f$ is any measurable function on $\bb T$, we can define the operator
of multiplication by $f$ from $L^2(\bb T)$ into the space of measurable functions.
For a subspace such as $H^2_{\alpha, \beta}$, the \textit{multipliers} are those
multiplication operators which map the space $H^2_{\alpha, \beta}$ into itself.
A standard argument using the Closed Graph Theorem shows that such operators
are always bounded.  We denote these multiplier operators by $M_f$.

\begin{prop} 
Let $0 \ne \alpha, \beta \in \bb C$ satisfy $|\alpha|^2+|\beta|^2=1.$ 
Then the algebra of multipliers of $H^2_{\alpha, \beta}$ is $H^{\infty}_1$,
and moreover $\|M_f\|= \|f\|_{\infty}$.
\end{prop}

\begin{proof} Let $f \cdot H^2_{\alpha, \beta} \subseteq H^2_{\alpha, \beta}.$ 
Since $z^2 \in H^2_{\alpha, \beta}$ it follows that $z^2f$ is analytic. 
So let $f(z) = \sum_{n \ge -2} f_n z^n$.
Observe that we must have $f(z)(\alpha + \beta z) \in H^2_{\alpha, \beta}.$ 
Examining the coefficient of $z^{-2}$ yields that $f_{-2}=0$; 
and similarly $f_{-1}=0.$
Equating the $0$th and first order terms yields
\[ f_0 \alpha + (f_0\beta+ f_1 \alpha) z = \lambda (\alpha + \beta z) \]
for some $\lambda$.
Hence   
$f_0= \lambda,$ and $f_0 \beta + f_1 \alpha = \lambda \beta = f_0 \beta,$ 
which yields $f_1=0.$ 
Thus every multiplier is in $H^{\infty}_1.$
But we noted earlier that every element of $H^{\infty}_1$ is a multiplier.

Finally, since $z^2H^2 \subseteq H^2_{\alpha, \beta} \subset H^2$,
\[
 \|f\|_\infty = \| M_f|_{z^2H^2}\| \le \| M_f|_{H^2_{\alpha, \beta}}\| 
 \le \| M_f|_{H^2}\| = \|f\|_\infty . \qedhere
\]
\end{proof}

We repeat the main theorem for the reader's convenience:

\renewcommand{\therepthm}{\ref{thm1}}
\begin{repthm}
Let $z_1,\dots,z_n$ be distinct points in $\bb D$, 
and let $w_1,\dots,w_n$ be complex numbers. 
Then there exists an analytic function $f$ on $\bb D$ 
with $\|f\|_{\infty} \le A$ and $f^{\prime}(0)=0$ such that
$f(z_i)=w_i$ for $i=1,\dots,n$ if and only if 
\[
 \Big[ (A^2-w_i \overline{w_j})K^{\alpha,\beta}(z_i,z_j) \Big] \ge 0 
 \qforal |\alpha|^2+|\beta|^2 =1.
\]
\end{repthm}

Now if $f$ is a multiplier in any reproducing kernel Hilbert space
with kernel $K(x,y)$ and $\|M_f\| \le A,$ 
then for any set of points, $x_1, \ldots, x_n$,
we will have that $\big[ (A^2 - f(x_i) \overline{f(x_j)}) K(x_i,x_j) \big]$
is positive semidefinite. See, for example, \cite[Theorem~5.2]{AM}.

Thus, if there exists $f \in H^{\infty}_1$ such that $f(z_i) = w_i$
and $\|f\|_{\infty} \le A$, then  
$f$ will be a multiplier of each $H^2_{\alpha, \beta}$ space 
with multiplier norm $\|M_f\| \le A$.
Hence 
\[ \Big[ (A^2 - w_i \overline{w_j}) K_{\alpha, \beta}(z_i,z_j) \Big] \ge 0 .\] 

This proves the necessity of the condition in Theorem \ref{thm1}.
The proof of sufficiency will follow from a factorization lemma 
similar to one used in \cite{Ab}.

Suppose that we are given a finite set $F= \{z_1,\dots,z_n \}$ of distinct points in $\bb D.$
Let $\cl I^F$ denote the ideal of functions in $H^\infty_1$ which
vanish on the set $F$.  
Write $B_F$ for the finite Blaschke product with simple zeroes at the points of $F$;  
i.e.\ $B_F = \prod_{z_i \in F} \phi_{z_i} = c_0 + c_1 z + \dots$.
Let $(a,b)$ be a norm one multiple of $(c_0,c_1)$.

When $0$ is not in $F$, it is evident that
\[ \cl I^F = H^\infty_1 \cap B_F H^\infty = z^2 B_F H^\infty + \bb C f_0 \]
where $f_0$ is the product of $B_F$ with that linear function which
multiplies it into $H^\infty_1$. 
It is easy to see that $f_0 = B_F (a-bz)$.
Therefore we conclude that 
\[ \cl I^F = B_F H^\infty_{a,-b} \]
where $H^\infty_{\alpha,\beta}$ denotes those $H^\infty$ functions
such that the vector $(h(0), h'(0))$ is a multiple of $(\alpha,\beta)$.
When $0$ belongs to $F$, write $F = F' \cup \{0\}$.  
Then $B_F = z B_{F'}$ and 
\[ \cl I^F = H^\infty_1 \cap B_F H^\infty = z^2 B_{F'} H^\infty = B_F H^\infty_{0,1} .\]
So the same result holds.

\begin{lemma} \label{L:pre-ann}
The pre-annihilator of $\cl I^F$ in $L^1(\bb T)$ is
$\cl I^F_\perp = \ol{z} \ol{B_F} H^1_{a,b}$.
\end{lemma}

\begin{proof}
Clearly $\cl I^F_\perp$ is contained in 
\[ (z^2B_F H^\infty)_\perp = \ol{z}^2 \ol{B_F} H^1_0 = \ol{z} \ol{B_F} H^1 ;\]
and consists of those functions in this set which are orthogonal to $f_0$.
(In the case of $0 \in F$, $f_0 = zB_F$, which is consistent.)
We may write $g \in \cl I^F_\perp$ as $g = \ol{z} \ol{B_F} g_0$ where $g_0 \in H^1$.
Then calculate
\begin{align*}
 0 &= \int f_0 g = \int B_F ( a - bz) \ol{z} \ol{B_F} g_0 \\
    &= a \int \ol{z} g_0 - b \int g_0 = a g_0'(0) - b g_0(0).
\end{align*}
Therefore $g_0$ belongs to $H^1_{a,b}$.
\end{proof}

Define $\cl N^F_{\alpha, \beta}$ to be the subspace of
$H^2_{\alpha, \beta}$ consisting of those functions which vanish on $F$.
Denote the orthogonal complement of $\cl N^F_{\alpha, \beta}$
in $H^2_{\alpha, \beta}$ by $\cl M^F_{\alpha, \beta}$.
This is $n$-dimensional, and evidently contains the kernel functions 
$k_{z_i}^{\alpha,  \beta}(z) = K_{\alpha, \beta}(z, z_i)$ for the points $z_i \in F$.
As they are linearly independent, these vectors span $\cl M^F_{\alpha, \beta}$.

\begin{lemma} \label{L:factor}
For each $g \in \cl I^F_\perp$, there are scalars $\alpha\ne 0$ and $\beta$ in $\bb C$
so that $g$ factors as $g= \ol{h} k$ where $k \in H^2_{\alpha,\beta}$,
$h \in L^2(\bb T)$ is orthogonal to $\cl N^F_{\alpha, \beta}$, and
$\|g\|_1= \|h\|_2 =\|k\|_2$.

Conversely, every function $g \in L^1(\bb T)$ which factors as 
$g= \ol{h} k$ where $k \in H^2_{\alpha,\beta}$ and
$h \in L^2(\bb T)$ is orthogonal to $\cl N^F_{\alpha, \beta}$
belongs to $\cl I^F_\perp$.
\end{lemma}

\begin{proof} 
Write $g = \ol{z} \ol{B_F} g_0$ as before and factor $g_0 = k h_0$
with 
\[ \|k\|_2 =\|h_0\|_2 = \|g_0\|_1 = \|g\|_1 ,\] 
where $k,h_0$ are in $H^2$ and $k$ is outer.
As $k = k_0 + k_1 z + \dots$ is outer, $k_0 \ne 0$.
Set $\alpha = k_0 (|k_0|^2+|k_1|^2)^{-1/2}$ and 
$\beta= k_1 (|k_0|^2+|k_1|^2)^{-1/2}$; 
so $\alpha \ne 0$ and $k$ lies in $H^2_{\alpha, \beta}$.
Write $h = z B_H \ol{h_0}$; so that $g = \ol{h} k$.
It remains to verify the orthogonality condition.

Observe that $H^\infty_1 k$ is dense in $H^2_{\alpha,\beta}$ because
outer functions are cyclic vectors for $H^\infty$ in $H^2$; whence
\[
 \ol{H^\infty_1 k} = \bb C k + \ol{ z^2 H^2 k} 
 = \bb C k + z^2 H^2 = H^2_{\alpha,\beta} .
\]
Also note that $\cl N^F_{\alpha, \beta} = \ol{\cl I^F H^2_{\alpha, \beta}}$.
Indeed, the right hand side vanishes on $F$, and so is contained in
$\cl N^F_{\alpha, \beta}$.  
However both evidently have codimension $n = |F|$; so they are equal.
Thus a dense subset of $\cl N^F_{\alpha, \beta}$ is given by $\cl I^F k$.

Therefore to check that $h$ is orthogonal to $\cl N^F_{\alpha, \beta}$, 
take an arbitrary element $f \in \cl I^F$ and calculate
\begin{align*}
 \inp{fk}{h} &= \inp{fk}{zB_F \ol{h_0}} = \int f(\ol{z}\ol{B_F} k h_0) = \int fg = 0 .
\end{align*}

Reversing this calculation shows that every such product $g=\ol{h}k$
belongs to $\cl I^F_\perp$.
\end{proof}

We return to the proof of the main result.

\begin{proof}[Proof of Theorem \ref{thm1}] 
It remains to be shown that if 
\[
 \Big[ (A^2 - w_i \overline{w_j})K_{\alpha, \beta}(z_i,z_j) \Big] \ge 0
 \qforal |\alpha|^2 + |\beta|^2 = 1 ,
\]
then there exists $\psi  \in H^{\infty}_1$ with $\|\psi\|_{\infty} \le A$ and
$\psi(z_i)=w_i$ for $i=1,\dots,n.$

Since the algebra of polynomials satisfying $p^{\prime}(0)=0$
separates points on the disk, we may choose a polynomial $p$ with
$p^{\prime}(0)=0$ and $p(z_i)=w_i$ for $i=1,\dots,n.$ 

For every $(\alpha, \beta)$, we have that 
$M_p^* k_{z_i}^{\alpha,\beta}= \overline{w_i} k_{z_i}^{\alpha,\beta}.$ 
Thus, $\cl M^F_{\alpha,\beta}$ is invariant under $M_p^*$.
As the functions $k_{z_i}^{\alpha,\beta}$ span $\cl M^F_{\alpha,\beta}$,
the positive semidefiniteness of the matrix above 
is equivalent to the condition 
\[ P_{\cl M^F_{\alpha,\beta}} (A^2 I - M_p M_p^*) |_{\cl M^F_{\alpha,\beta}} \ge 0;\]
which is equivalent to $\| M_p^*|_{\cl M^F_{\alpha,\beta}} \| \le A$.

Let $P_{\alpha,\beta}: L^2(\bb T) \to H^2_{\alpha,\beta}$
denote the orthogonal projection.

Define a linear functional  $\Phi$ on $\cl I^F_\perp$ by $\Phi(f) = \int p f$.
We claim that $\Phi$ has norm at most $A$. 
To see this, take any $g \in \cl I^F_\perp$ and factor it as $g = \ol{h}k$
as in Lemma~\ref{L:factor}. 
Since $h$ is orthogonal to $\cl N^F_{\alpha,\beta}$,
we see that $P_{\alpha,\beta}(h)=:\tilde h$ lies in $\cl M^F_{\alpha,\beta}$.
We compute
\begin{align*}
 \Phi(g) &= \int pk\ol{h} = \inp{pk}{h} = \inp{P_{\alpha,\beta}pk}{h} \\
 &= \inp{pk}{P_{\alpha,\beta} h} = \langle k,M_p^* \tilde h \rangle .
\end{align*}
Therefore 
\[
 |\Phi(g)| \le \|M_p^*|_{\cl M^F_{\alpha,\beta}} \|\,  \|\tilde h\|_2 \, \|k\|_2 \le A  \|g\|_1 .
\]

Thus, by the Hahn-Banach Theorem, we may extend $\Phi$ to a linear functional on
$L^1$ of norm at most $A$. 
Since $L^{\infty}$ is the dual of $L^1$, this means that there is a function 
$f \in L^{\infty}$ with $\|f \|_{\infty} \le A$ such that
$\Phi(g) = \int f g = \int pg$ for every $g \in \cl I^F_{\alpha,\beta}.$
Therefore $f - p$ belongs to $(\cl I^F_\perp)^\perp = \cl I^F$.
In particular, $f \in H^\infty_1$ and $f-p$ vanishes on $F$.
So $f(z_i) = p(z_i) = w_i$ as desired.
\end{proof}

The proof of Theorem \ref{thm3} is easier.

\renewcommand{\therepthm}{\ref{thm3}}
\begin{repthm}
Let $z_1,\dots,z_n$ be distinct points in $\bb D$, 
and let $w_1,\dots,w_n$ be in $\bb D$. Then there exists an
analytic function $f \in H^\infty_1$ with $\|f\|_{\infty} \le 1$ 
such that $f(z_i)=w_i$ for $i=1,\dots,n$ if and only if
there exists $\lambda \in \bb D$ so that
\[ \left[ \frac{z_i^2 \overline{z_j^2} - \phi_{\lambda}(w_i)
  \overline{\phi_{\lambda}(w_j)}}{1-z_i \overline{z_j}} \right] \ge 0.\]
\end{repthm}

\begin{proof}
We first suppose that $0 \not\in F$.
Assume that $f \in H^\infty_1$ exists satisfying $\|f\|_{\infty} \le 1$ 
and $f(z_i)=w_i$ for $i=1,\dots,n$. 
Set $\lambda = f(0)$, which lies in $\bb D$; and let $g(z)
= \phi_{\lambda}(f(z)).$ 
Then $g$ belongs to $H^\infty$ and $g'(0)=\phi_\lambda'(f(0)) f'(0) = 0$.
So $g \in H^\infty_1$.
Evidently $\|g\|_{\infty} \le 1$, $g(z_i) = \phi_{\lambda}(w_i)$ for $i=1,\dots,n$ and
$g(0)=\phi_\lambda(\lambda)=0.$ 
Hence, $g(z)=z^2h(z)$ with $\|h\|_{\infty} \le 1.$ 

Therefore by the Nevanlinna--Pick Theorem applied to $h$ and the set $F$, 
\[ 
 \left[ \frac{1 - h(z_i)\ol{h(z_j)}} {1-z_i\ol{z_j}} \right] \ge 0.
\]
Consequently, $\left[ \frac{z_i^2 \ol{z_j^2} - g(z_i) \ol{g(z_j)}}{1-z_i \ol{z_j}} \right]$
equals
\[
 \begin{bmatrix} z_1^2 & 0 & \dots & 0\\0 & z_2^2 & \dots & 0\\
 \vdots & \vdots & \ddots & \vdots\\ 0 & 0 & \dots & z_n^2 \end{bmatrix}
 \left[ \frac{1 - h(z_i)\ol{h(z_j)}} {1-z_i\ol{z_j}} \right] 
 \begin{bmatrix} \ol{z_1}^2 & 0 & \dots & 0\\0 & \ol{z_2}^2 & \dots & 0\\ 
 \vdots & \vdots & \ddots & \vdots\\ 0 & 0 & \dots & \ol{z_n}^2 \end{bmatrix} ;
\]
and thus is positive semidefinite. 

On the other hand, if, say $z_1=0$, then $\lambda=f(0)=f(z_1)=w_1$ and
the first row and column of the matrix above is zero; and the same reasoning
applies to the remaining entries of the matrix.

Conversely, we again consider the case where $0 \not\in F$, and suppose that
$\lambda$ is given which provides a positive semidefinite matrix.
By reversing the calculation, we see that the matrix, 
\[
 \left[\frac
 {1- z_i^{-2}\phi_{\lambda}(w_i) \ol{z_j^{-2}\phi_{\lambda}(w_j)}}
 {1-z_i \ol{z_j}} \right]
\]
is positive semidefinite.
Hence, by the Nevanlinna--Pick Theorem, there exists an analytic function $h$ on
$\bb D$ with $\|h\|_{\infty} \le 1$ so that $h(z_i) = z_i^{-2}\phi_{\lambda}(w_i).$ 
Set $f(z) = \phi_{-\lambda}(z^2h(z)).$ 
Reversing the calculations of the first paragraph shows that
$f$ is the desired interpolant. 

The case where $0 \in F$ is handled similarly.
\end{proof}

\section{Distance Formulae} \label{S:dist}

Donald Sarason introduced new operator theoretic methods for interpolation
problems in his seminal paper  \cite{Sa}.
In particular, he made critical use of the fact that the Nevanlinna--Pick 
interpolation problem is equivalent to a distance estimate.
In our context, we are searching for a function $f \in H^\infty_1$
of minimal norm satisfying $f(z_i)=w_i$ for $z_i \in F$, $i=1,\dots,n$.
Letting $p$ be a polynomial in $H^\infty_1$ satisfying the interpolation data,
it is easy to see that the optimal norm is precisely $\dist(p,\cl I^F)$.
The infimum is attained because $\cl I^F$ is weak-$*$ closed.

The main theorem can be re-interpreted as the following.

\begin{thm}\label{T:dist ideal}
Let $F$ be a finite subset of $\bb D$.  For any $f \in H^\infty_1$,
\[
 \dist(f, \cl I^F) =
 \sup_{|\alpha|^2+|\beta|^2=1} \|M_f^*|_{\cl M^F_{\alpha,\beta}}\| .
\]
\end{thm}

\begin{proof}
We compute $\dist(f, \cl I^F)$ using duality.
By Lemmas~\ref{L:pre-ann} and \ref{L:factor}, we have
\begin{align*}
 \dist(f, \cl I^F) &= 
 \sup_{\substack{g \in \cl I^F_\perp\\ \|g\|_1 = 1}}
 \Big| \int fg\ \Big| =
 \sup_{\substack{|\alpha|^2+|\beta|^2=1\\
 k \in H^2_{\alpha,\beta},\ \|k\|_2\le 1\\ 
 h \perp \cl N^F_{\alpha,\beta},\ \|h\|_2\le 1}}
 \Big| \int fk\ol{h}\ \Big| \\ &= \sup \big| \inp{k}{M_f^* P_{\alpha,\beta}h} \big|
 = \sup_{|\alpha|^2+|\beta|^2=1} \|M_f^*|_{\cl M^F_{\alpha,\beta}}\| .
 \qedhere
\end{align*}
\end{proof}

We recover Theorem~\ref{thm1} by observing that $\|M_f^*|_{\cl M^F_{\alpha,\beta}}\| \le A$
if and only if 
\[
 \Big[ (A^2-w_i \overline{w_j})K^{\alpha,\beta}(z_i,z_j) \Big] \ge 0 .
\]

Something special happens when $F$ contains $0$.
In this case, $\cl I^F = zB_F H^\infty$.
This is an ideal of $H^\infty$ as well, so  the distance to it can be 
computed using standard methods.  
Let $\Gamma_h$ denote the Hankel operator $P^\perp M_f |_{H^2}$
where $P$ is the projection onto $H^2$.
Using Nehari's Theorem,
\begin{align*}
 \dist(f,\cl I^F) &= \dist(\ol{z}\ol{B_F}f, H^\infty) = \| \Gamma_{\ol{z}\ol{B_F}f} \| \\
 &= \| P^\perp M_{zB_F}^* M_fP\| = \| P^\perp M_{zB_F}^* PM_fP\| \\
 &= \| P M_f^* PM_{zB_F} P^\perp\| = \|PM_f^*|_{H^2 \ominus zB_F H^2}\| .
\end{align*}
Therefore, we obtain:

\begin{cor} \label{C:dist0}
Let $F=\{z_1 = 0,z_2,\dots,z_{n}\}$ be a finite subset of $\bb D$ containing $0$.
Given data $w_i \in \bb C$ for $i=1,\dots,n$, the following are equivalent:
\begin{enumerate}
\item there is a function $f$ in $H^\infty_1$
  such that $f(z_i) = w_i$ for $i=1,\dots,n$ and $\|f\|_\infty \le A$.\vspace{.5ex}

\item $\dist(f,\cl I^F) \le A$.
i.e.\ $\dist(f,zB_F H^\infty) \le A$.
\vspace{.5ex}

\item $\| PM_f^*|_{H^2 \ominus zB_F H^2}\| \le A$.
\vspace{.5ex}

\item {\renewcommand{\arraycolsep}{.5pt}
  $\left[ \begin{array}{cc|ccc}
  A^2 \!-\! |w_1|^2&0&A^2 \!-\! w_1\ol{w_2}&\dots&A^2 \!-\! w_1\ol{w_n}\\
  0&A^2 \!-\! |w_1|^2&(A^2 \!-\! w_1\ol{w_2})\ol{z_2}&\dots&(A^2 \!-\! w_1\ol{w_n})\ol{z_n}\\[.5ex]
  \hline \phantom{\Big|} 
  A^2 \!-\! w_2 \ol{w_1}\ &(\!A^2 \!-\! w_2 \ol{w_1})z_2 &
  \frac{A^2-w_2 \ol{w_2}}{1-z_2\ol{z_2}} &\dots& \frac{A^2-w_2 \ol{w_n}}{1-z_2\ol{z_n}}\\
  \vdots&\vdots&\vdots&\frac{A^2-w_i \ol{w_j}}{1-z_i\ol{z_j}} &\vdots\\
  A^2 \!-\! w_n \ol{w_1}\ &(\!A^2 \!-\! w_n \ol{w_1})z_n&
  \frac{A^2-w_n \ol{w_2}}{1-z_n\ol{z_2}} &\dots& \frac{A^2-w_n \ol{w_n}}{1-z_n\ol{z_n}}
  \end{array}\right]$\\[1ex]
  } 
  is positive semidefinite.
\end{enumerate}
\end{cor} 

\begin{proof}
We have already shown the equivalence of (1), (2) and (3).
To obtain the equivalence of (3) and (4), we use the basis
$1,z,k_{z_2},\dots,k_{z_n}$ for $H^2 \ominus zB_F H^2$,
where $k_{z_i}(z) = \frac1{1-\ol{z_i}z}$ are the reproducing kernel
functions for the Hardy space $H^2$.
Observe that since $f\in H^\infty_1$, we have $M_f^*z = \ol{w_0} z$.
The positivity of $P_{H^2 \ominus zB_F H^2} (A^2 - M_f M_f^*) P_{H^2 \ominus zB_F H^2}$
is equivalent to the positivity of the $(n+1)\times(n+1)$ matrix
$\Big[ \big\langle (A^2 - M_f M_f^*) k_j,k_i \big\rangle \Big]$ where 
$k_0=1$, $k_1=z$ and $k_i = k_{z_i}$ for $i=2,\dots,n$.
A simple computation shows that this is the matrix in the semidefinite condition (4).
\end{proof}

It is interesting to reconcile this corollary with Theorem~\ref{T:dist ideal}.
Notice that when $0$ is in $F$, that 
\[ \cl N^F_{\alpha,\beta} = H^2_{\alpha,\beta} \cap B_F H^2 = zB_F H^2 \]
independent of $(\alpha,\beta)$.  Therefore
\[ \cl M^F_{\alpha,\beta} = ( z^2 H^2 \ominus zB_F H^2) \oplus \bb C (\alpha + \beta z) .\]
The special thing that occurs here is that the subspaces $\cl M_{\alpha,\beta}$
have codimension 1 in the space $H^2 \ominus zB_F H^2$.
Moreover, any vector $x$ in $H^2 \ominus zB_F H^2$ has the form
$x = c(\alpha + \beta z) + z^2h$; and therefore lies in one of the
subspaces $\cl M^F_{\alpha,\beta}$.
Consequently, if $M_f^*|_{H^2 \ominus zB_F H^2}$ achieves its norm
at a vector $x$, and $x \in \cl M^F_{\alpha,\beta}$, then the same norm
is achieved on the restriction $M_f^*|_{\cl M^F_{\alpha,\beta}}$.
That is, when $0 \in F$,
\[
 \sup_{|\alpha|^2+|\beta|^2=1} \|M_f^*|_{\cl M^F_{\alpha,\beta}}\|
 = \| PM_f^*|_{H^2 \ominus zB_F H^2}\| .
\]

In the case in which $0 \not\in F$, $\cl N^F_{\alpha,\beta}$ contains
$z^2B_F H^2$ as a codimension 1 subspace.  The subspaces
$\cl M^F_{\alpha,\beta}$ are codimension 2 in their joint span, 
not codimension 1.  So the argument above does not apply.

The implications of this are explored in the discussion of C*-envelopes.

In the same vein as Theorem~\ref{T:dist ideal}, 
we can obtain an analogue of Nehari's Theorem:
if $f \in L^\infty$, then $\dist(f,H^\infty) = \| P^\perp M_f P\| = \| \Gamma_f\|$.

\begin{thm} \label{distform}
If $f\in L^\infty$, then 
\[
 \dist(f, H^\infty_1) =
 \sup_{|\alpha|^2+|\beta|^2=1}
 \|(I-P_{\alpha,\beta})M_f P_{\alpha,\beta}\|.
\] 
\end{thm}

\begin{proof}
We use a similar duality argument. 
The pre-annihilator of $H^\infty_1$ in $L^1$ is the closed span of 
$H^1_0$ and $\ol{z}$; and this is equal to $\ol{z} H^1_1$.
This is just Lemma~\ref{L:pre-ann} for $F=\emptyset$.
Then Lemma~\ref{L:factor} shows that every $g$ in $(H^\infty_1)_\perp$ factors
as $g = \ol{h}k$ where $\|k\|_2=\|h\|_2=\|g\|_1$, $k$ lies in some $H^2_{\alpha,\beta}$
and $h$ is orthogonal to $\cl N^\emptyset_{\alpha,\beta} = H^2_{\alpha,\beta}$.
Conversely, every product of this form lies in $(H^\infty_1)_\perp$.
Therefore
\begin{align*}
 \dist(f, H^\infty_1) &= \sup_{g \in (H^\infty_1)_\perp,\, \|g\|_1=1} \Big| \int fg \ \Big|
 = \sup_{\substack{|\alpha|^2+|\beta|^2=1\\
 k \in H^2_{\alpha,\beta},\ \|k\|_2\le 1\\ 
 h \perp \cl H^2_{\alpha,\beta},\ \|h\|_2\le 1}} 
 \big| \inp{fk}{h} \big| 
 \\ &= \sup_{|\alpha|^2+|\beta|^2=1}
 \|(I-P_{\alpha,\beta})M_f P_{\alpha,\beta}\|.
 \qedhere
\end{align*}
\end{proof}

\section{Matrix-Valued Interpolation and C*-envelopes} \label{S:Cenv}

In the classical Nevanlinna--Pick problem, one can consider matrix
valued interpolation.  That is, one specifies points $z_1,\dots,z_n$
in the unit disk and $k\times k$ matrices $W_1,\dots,W_n$ and asks
for the optimal norm $\|f\|_\infty$ of a bounded analytic function $f$
from $\bb D$ into $\ff M_k$ satisfying $f(z_i) = W_i$.  
The norm of a function in $\ff M_k(H^\infty)$ is defined as the supremum 
over $\bb D$ of the operator norm of the matrix $f(z)$.
It turns out that the same result holds,
namely that there is such a function with $\|f\|\le A$ if and only if
the matrix $\Big[ \frac{A^2 I_k - W_i W_j^*}{1-z_i\ol{z_j}} \Big]$ is
positive semidefinite.  

The same matrix interpolation problem can be formulated for 
the space of multipliers of any reproducing kernel Hilbert space.  
When the positivity of the matrix 
$\big[ (A^2 I_k - W_i W_j^*)K(z_i,z_j) \big]$ is equivalent to interpolation,
the kernel is called a {\em complete} Nevanlinna--Pick kernel.
Such kernels have been characterized by Agler and McCarthy \cite{AM3}.
However, recognition of such kernels is generally not straightforward.

In our context of $H^\infty_1$, the analogous problem is to ask whether 
the family of conditions
\[ 
 \big[(A^2 I_k - W_iW_j^*) K^{\alpha,\beta}(z_i,z_j) \big] \ge 0
 \qforal |\alpha|^2+|\beta|^2=1 
\]
is equivalent to matrix interpolation of the data by a function $f\in \ff M_k(H^\infty_1)$
with $\|f\| \le A$.

When the set consists of only two-points, the answer to such questions is always affirmative, for reasons that we shall discuss in section 7. However, for three or more points the problem is more difficult.

Recall that every unital operator algebra $\cl A$ imbeds completely
isometrically into some C*-algebra.  Moreover, among such C*-algebras
which are generated by the range, there is a unique smallest one
known as the C*-envelope, $\ce(\cl A)$, in the sense that:
if $j$ is a completely isometric isomorphism of $\cl A$ into a 
C*-algebra $\ff A = \ca(j(\cl A))$, then there is an ideal $\ff I$ of $\ff A$
so that the quotient map $q$ by $\ff I$ is a complete isometry 
on $j(\cl A)$ and $\ff A/\ff I$ is $*$-isomorphic to $\ce(\cl A)$.
See \cite{Pau} for the necessary background.

One way to determine the complexity of the matrix interpolation
problem is to compute the C*-envelope of the associated quotient algebra. 
This connection between interpolation and the C*-envelope has
been studied and is further discussed in \cite{McP} and \cite{So}.
In this section, we use the computation of a C*-envelope 
to show that for certain subsets of $\bb D$,
this matrix-valued analogue fails.
To be precise, we will prove the following result.

\begin{thm} \label{T:fail matrix} 
There exists a set $F=\{z_1, z_2, z_3\}$ of three distinct non-zero
points in $\bb D,$ an integer $k$ and $k \times k$ matrices $W_1, W_2, W_3$ 
such that the $3 \times 3$ block matrix $\big[(I_k - W_iW_j^*) K^{\alpha,\beta}(z_i,z_j) \big]$
is positive semidefinite for all $|\alpha|^2+ |\beta|^2 =1$, but there
does not exist a function $f \in \ff M_k(H^\infty_1)$ with $\|f\|_\infty \le 1$ 
such that $f(z_i)=W_i$ for $i=1,2,3$.
\end{thm}

Our proof is indirect and, in particular, we are currently unable to
explicitly exhibit a particular set of three points and the three matrices
$W_1, W_2, W_3,$ which the above theorem asserts exist.

Let $F = \{z_1,\dots,z_n \}$ be a finite subset of $\bb D$
and let $W_1,\dots,W_n$ belong to $\ff M_k$.
Consider the problem of finding the optimal norm of a function
$f \in \ff M_k(H^\infty_1)$ satisfying $f(z_i) = W_i$.
Such functions always exist, even amongst polynomials.
Let $p$ be an arbitrary choice of an interpolant.
As in the previous section, the matrix interpolation problem is
equivalent to distance estimate.
The optimal norm is $A := \dist(p,\ff M_k(\cl I^F))$. Because of this equivalence between distance and interpolation, we can re-interpret the above theorem in terms of distance formulae.

\begin{cor} \label{fail distance} There exists a set $F=\{z_1, z_2, z_3 \}$ of three distinct non-zero points in $\bb D,$ an integer $k$ and a $k \times k$ matrix-valued function $p \in \ff M_k(H^\infty_1)$ such that

\[
 \dist(p, \ff M_k(\cl I^F)) \ne
 \sup_{|\alpha|^2+|\beta|^2=1} \|M_p^*|_{\bb C^k \otimes \cl M^F_{\alpha,\beta}}\| .
\]
\end{cor}

It is likely that the matrix-valued analogue of \ref{distform} is also false although we have not shown that here.

To prove the above theorem, we need to first
reinterpret the results of the previous section.
It is now important to use each kernel exactly once, so we use the 
parameterization by the projective complex 2-sphere $PS^2$.
But for convenience of notation, we still use $(\alpha,\beta)$
to represent a point in $PS^2$.

Define a map $\Phi_F$ from $H^\infty_1$ into $\rm C(PS^2,\ff M_n)$ by
\[
 \Phi_F(f)(\alpha,\beta) = 
 P_{\cl M^F_{\alpha,\beta}}M_f|_{\cl M^F_{\alpha,\beta}} .
\]
We use the compression of $M_f$ rather than the restriction of
$M_f^*$ so that our map is linear.
The image is a continuous function because the map taking 
$(\alpha,\beta)$ to the projection $P_{\cl M^F_{\alpha,\beta}}$ is continuous.
Clearly $\ker \Phi_F = \cl I^F$.  
Therefore it induces a map $\tilde\Phi_F$ from $H^\infty_1/\cl I^F$
into $\rm C(PS^2,\ff M_n)$.
Theorem~\ref{T:dist ideal} says that $\tilde\Phi_F$ is isometric.

For convenience, we write $\cl K := H^2\ominus zB_F H^2$.
When $F$ contains $0$, Corollary~\ref{C:dist0} provides a different map
$\Psi_F$ from $H^\infty_1$ into $\ff M_{n+1}$ given by
\[ \Psi_F(f) = P_{\cl K} M_f|_{\cl K} .\]
Again this factors through the quotient by $\cl I^F$, and yields an isometric
map $\tilde\Psi_F$ from $H^\infty_1/\cl I^F$ into $\ff M_{n+1}$.

However, in this case, more is true.
The map $\Psi_F$ extends naturally to $H^\infty$,
and we keep the same name for it.
Since $0 \in F$,  $\cl I^F = zB_F H^\infty$.
Therefore the injection of $H^\infty_1/\cl I^F$ into 
$H^\infty/zB_F H^\infty$ is completely isometric.
The classical matrix Nevanlinna--Pick interpolation result
is equivalent to saying that the map $\tilde\Psi_F$ is a
\textit{complete} isometry from $H^\infty/zB_F H^\infty$
into $\ff M_{n+1}$.
A fortiori, the map $\tilde\Psi_F$ restricted to $H^\infty_1/\cl I^F$
is a complete isometry.

This will enable us to compute the C*-envelope of 
$H^\infty_1/\cl I^F$ in this case.
Indeed, the algebra $\Psi_F(H^\infty)$ is known to generate
all of $\ff M_{n+1}$ as a C*-algebra.
Since $\ff M_{n+1}$ is simple, it is the C*-envelope
of $H^\infty/zB_F H^\infty$.  
Usually this is the case for $H^\infty_1$ as well, with $n=2$
being an exception.

\begin{thm} \label{T:Cenv0}
Let $F$ be a set of $n$ distinct points in $\bb D$ containing $0$.
If $n\ge 3$, then $\ce(H^\infty_1/\cl I^F) = \ff M_{n+1}$. 
If $n=2$, then $\ce(H^\infty_1/\cl I^F) = \ff M_2$.
\end{thm}

\begin{proof}
We first need a useful representation of $\Psi_F(f)$.
As in the proof of Corollary~\ref{C:dist0}, we make use of the basis
$1,z,k_{z_2},\dots,k_{z_n}$ for $\cl K$, and the fact that 
$\Psi_F(f)^* 1 = \ol{w_1}$, $\Psi_F(f)^* z = \ol{w_1}z$ and
$\Psi_F(f)^* k_{z_i} = \ol{w_i} k_{z_i}$ for $i=2,\dots,n$,
where $f(z_i)=w_i$.  
Thus $\Psi_F(f)^*$ is diagonal with respect to this non-orthogonal basis.  

Let $D_f = \diag \big( f(z_1), f(z_1), f(z_2),\dots,f(z_n) \big)$ be the
diagonal $n\!+\!1 \times n\!+\!1$ matrix in $\ff M_{n+1}$ with the first 
eigenvalue repeated a second time.
It will be convenient to write the standard basis of $\ff M_{n+1}$ as
$e_0,\dots,e_n$.
Consider $V = \big[1,z,k_{z_2},\dots,k_{z_n} \big]^*$ as a map from
$\cl K$ into $\bb C^{n+1}$. 
Then $\Psi(f)^* = V^{-1}D_f^*V$.
Using polar decomposition, we may replace $V$ with the map 
$Q^{1/2}$ where
\begin{align*}
Q = VV^*  &=
\begin{bmatrix}
 \inp{1}{1} & \inp{1}{z} & \inp{1}{k_{z_2}} &\ldots & \inp{1}{k_{z_n}} \\ 
 \inp{z}{1} & \inp{z}{z} & \inp{z}{k_{z_2}} &\ldots & \inp{z}{k_{z_n}} \\
 \inp{k_{z_2}}{1} & \inp{k_{z_2}}{z} \\
 \vdots & \vdots & & \big[ \inp{k_{z_j}}{k_{z_i}} \big]_{i,j\geq 2}\\ 
 \inp{k_{z_n}}{1} & \inp{k_{z_n}}{z}
\end{bmatrix} \\ &=
\begin{bmatrix}
 1 & 0 & 1 &\ldots & 1 \\
 0 & 1 & z_2 &\ldots & z_n \\
 1 & \overline{z_2} \\
 \vdots & \vdots & & \left[\dfrac{1}{1-z_i\overline{z_j}}\right]_{i,j\geq 2}\\ 
 1 & \overline{z_n}
\end{bmatrix}
\end{align*}
Hence $\Psi(f)^* \simeq Q^{-1/2} D_f^* Q^{1/2}$;
and so $\Psi(f) \simeq Q^{1/2} D_f Q^{-1/2}$.
Therefore $\Psi$ is unitarily equivalent to the map
$\pi(f) = Q^{1/2} D_f Q^{-1/2}$.

Observe that $H^\infty_1/\cl I^F$ is generated by the $n$ commuting 
idempotents which are cosets of functions $f_j$ satisfying $f_j(z_i) = \delta_{ij}$
for $1 \le i,j \le n$.  
Write $E_{ij}$, $0 \le i,j \le n$, for the matrix units in $\ff M_{n+1}$.
The functions $f_j$ are mapped to $\pi(f_1) = Q^{1/2}(E_{0,0}+E_{1,1})Q^{-1/2}$
and $\pi(f_j) = Q^{1/2}(E_{j,j})Q^{-1/2}$ for $j=2,\dots,n$. 

Let $\ff A = \ca(\pi(H^\infty_1))$.
Then $\ff A$ must contain the operators $\pi(f_j)^* \pi(f_j)$.
Observe that 
\begin{align*}
 \pi(f_1)^* \pi(f_1) &= 
 Q^{-1/2}(E_{0,0}+E_{1,1})Q(E_{0,0}+E_{1,1})Q^{-1/2} \\&=
 Q^{-1/2}(E_{0,0}+E_{1,1})Q^{-1/2}
\end{align*}
and for $j=2,\dots,n$, 
\[ \pi(f_j)^* \pi(f_j) = Q^{-1/2} E_{j,j} Q E_{j,j} Q^{-1/2} = q_{jj} Q^{-1/2} E_{j,j} Q^{-1/2} \]
where $q_{jj} = (1-|z_j|^2)^{-1}$ is the $j,j$ entry of $Q$.
In particular, $\ff A$ contains
\begin{align*}
\sum_{j=1}^n q_{jj}^{-1} \pi(f_j)^* \pi(f_j) &= 
 Q^{-1/2} (E_{0,0} \!+\! E_{1,1})Q^{-1/2} 
+\sum_{j=2}^n  Q^{-1/2}E_{j,j}Q^{-1/2} \\
&= Q^{-1/2} I  Q^{-1/2} =  Q^{-1} .
\end{align*}
Hence $Q$, $Q^{1/2}$ and $Q^{-1/2}$ all belong to $\ff A$. 

Therefore $E_{0,0}+E_{1,1}$ and $E_{j,j}$ for $j=2,\dots,n$ belong to $\ff A$.
In addition, $E_{ii}QE_{jj} = q_{ij}E_{ij}$ is a non-zero multiple of $E_{ij}$
in $\ff A$ for $2 \le i,j \le n$.
And similarly, $\ff A$ contains 
\[ (E_{0,0}+E_{1,1}) Q E_{jj} = E_{0j} + z_j E_{1j} \qfor j = 2,\dots,n .\]

As long as $n \ge 3$, we obtain $E_{0j}$ and $E_{1j}$ in $\ff A$.
Therefore $\ff A$ is all of $\ff M_{n+1}$.
As this is a simple C*-algebra, it must be the C*-envelope.

Now consider $n=2$.
In this case $\ff A$ is generated by
\[ 
 Q = \begin{bmatrix} 1 & 0 & 1 \\ 0 & 1 & z_2 \\ \ol{z_2} & 1 & \frac1{1-|z_2|^2} \end{bmatrix},\
 D_1 = \begin{bmatrix} 1 & 0 & 0 \\ 0 & 1 & 0 \\ 0 & 0 & 0 \end{bmatrix} \AND
 D_2 = \begin{bmatrix} 0 & 0 & 0 \\ 0 & 0 & 0 \\ 0 & 0 & 1 \end{bmatrix} .
\]
Let
\[ 
C = (E_{0,0}+E_{1,1}) Q E_{22} = 
\begin{bmatrix} 0 & 0 & 1 \\  0 & 0 & z_2 \\  0 & 0 & 0 \end{bmatrix}.
\] 
Then $(1 + |z_2|^2)^{-1/2} C$ is a rank 1 partial isometry, and 
$\spn\{CC^*,C,C^*,D_2\}$ is a copy of $\ff M_2$ on 
the subspace $\cl N :=\spn\{(1,z_2,0), (0,0,1)\}$; 
while $E_0:= D_1 - (1 + |z_2|^2)^{-1}CC^*$ spans a copy of $\bb C$
on the complement $\spn\{(\ol{z_2},-1,0)\}$.
Since $Q = D_1 + D_2 + C + C^*$, it is evident that these two subalgebras
generate all of $\ff A$.
So $\ff A \simeq \ff M_2 \oplus \ff M_1$.

To see that the C*-envelope is $\ff M_2$ in this case, it suffices
to show that the quotient map onto the $\ff M_2$ summand is 
completely isometric on $H^\infty_1/\cl I^F$.
Let $E_1 = D_1 - E_0$ and $E_2 = D_2$, so that $E_0$, $E_1$ and $E_2$
are diagonal matrix units compatible with the decomposition of $\ff A$.
Observe that $Q = E_0 \oplus Q_1$ where $Q_1 = (E_2+E_3) Q (E_2+E_3)|_{\cl N}$.
Hence $Q^{1/2} = E_1 \oplus Q_1^{1/2}$.

A typical element of $\ff M_k(H^\infty_1/\cl I^F)$ has the form 
$X = f_1\otimes A_1 + f_2\otimes A_2$ for matrices $A_1$ and $A_2$ in $\ff M_k$.
Using the structure above, we calculate
\begin{align*}
 \pi(X) &= 
 (Q^{-1/2}\otimes I_k) ((E_0+E_1)\otimes A_1 + E_2\otimes A_2 ) (Q^{1/2}\otimes I_k) \\
 &= (E_0\otimes A_1) \oplus 
 (Q_1^{-1/2}\otimes I_k) (E_1\otimes A_1 + E_2\otimes A_2  ) (Q^{1/2}\otimes I_k) \\
 &= (E_0\otimes A_1) \oplus (P_1\otimes A_1 + P_2\otimes A_2 )
\end{align*}
where $P_1=P_1^2$ is idempotent and $P_2 = I - P_1$.
In order to show that the second summand always dominates the first
in norm, it suffices to show that
\[ \|A_1\| \le \| P_1\otimes A_1 + P_2\otimes A_2 \| .\]
To see this, take a unit vector $e$ in the range of $P_1$ and
a unit vector $x$ such that $\|A_1x\| = \|A_1\|$.  Then
\[
 \| P_1\otimes A_1 + P_2\otimes A_2 \| \ge
 \| (P_1\otimes A_1 + P_2\otimes A_2) e\otimes x \| =
 \| e \otimes A_1x \| = \|A_1\| .
\]
Therefore the quotient of $\ca(\pi(H^\infty_1/\cl I^F))$ onto $\ff M_2$
is completely isometric on $H^\infty_1/\cl I^F$; and hence $\ff M_2$
is the C*-envelope.
\end{proof}

\begin{cor} \label{C:not c.isom}
Let $F$ be a set of $n\ge3$ distinct points in $\bb D$ containing $0$. 
Then the isometric homomorphism 
$\tilde\Phi_F: H^\infty_1/\cl I^F \to \rm C(PS^2, \ff M_n)$ is not completely isometric.
\end{cor}

\begin{proof}
If $\tilde\Phi_F$ were a complete isometry, then there would 
be a *-homo\-mor\-phism from the C*-subalgebra $\ff A$ of 
$\rm C(PS^2, \ff M_n)$  generated by the range of 
$\tilde\Phi_F$ onto $C^*_e(H^\infty_1/\cl I^F) = \ff M_{n+1}.$ 
However, as $\ff A$ is a subalgebra of $\rm C(PS^2, \ff M_n)$,
every irreducible representation of $\ff A$ is of dimension at most $n$. 
This contradiction leads to the conclusion that $\tilde\Phi_F$ 
is not a complete isometry.
\end{proof}

\begin{cor} \label{C:interp0}
Let $z_1=0,z_2,\ldots,z_n$ be $n$ distinct points in the disk and let $W_1,\ldots,W_n\in \ff M_k$ be $k\times k$ matrices. There exists a function $f\in \ff M_k(H^\infty_1)$ with $\norm{f}\leq 1$ and $f(z_j)=W_j$ if and only if the matrix

\renewcommand{\arraycolsep}{.5pt}
  $$\left[ \begin{array}{cc|ccc}
  I \!-\! W_1W_1^\ast&0&I \!-\! W_1W_2^\ast&\dots&I \!-\! W_1W_n^\ast\\
  0&I \!-\! W_1W_1^\ast&(I \!-\! W_1W_2^\ast)\ol{z_2}&\dots&(I \!-\! W_1W_n^\ast)\ol{z_n}\\[.5ex]
  \hline \phantom{\Big|} 
  I \!-\! W_2W_1^\ast\ &(\!I \!-\! W_2W_1^\ast)z_2 &
  \frac{I-W_2W_2^\ast}{1-z_2\ol{z_2}} &\dots& \frac{I-W_2W_n^\ast}{1-z_2\ol{z_n}}\\
  \vdots&\vdots&\vdots&\frac{I-W_iW_j^\ast}{1-z_i\ol{z_j}} &\vdots\\
  I \!-\! W_nW_1^\ast\ &(\!I \!-\! W_nW_1^\ast)z_n&
  \frac{I-W_n W_2^\ast}{1-z_n\ol{z_2}} &\dots& \frac{I-W_nW_n^\ast}{1-z_n\ol{z_n}}
  \end{array}\right]$$
is positive semidefinite.
\end{cor}

Now we use this result to deduce that matrix interpolation must fail
for certain sets that do not contain $0$.

Let $F = \{z_1,\dots,z_n\}$ be a set of distinct non-zero points in $\bb D$.
As in the proof of the previous theorem, observe that $\Phi_F(f)(\alpha,\beta)$
is similar to the diagonal operator $D_f := \diag\big( f(z_1),\dots, f(z_n) \big)$ with
respect to the non-orthogonal basis $k_{z_1},\dots,k_{z_n}$.
The operator $V_{\alpha,\beta} = \big[ k_{z_1},\dots,k_{z_n} \big]^*$ from
$\cl M^F_{\alpha,\beta}$ to $\bb C^n$ implements the similarity via
$\Phi_F(f)(\alpha,\beta) = V_{\alpha,\beta}^{-1} D_f V_{\alpha,\beta}$.
As before, we define
\[
 Q_{\alpha,\beta} = V_{\alpha,\beta} V_{\alpha,\beta}^* 
 = \Big[ \ip{k_{z_j}, k_{z_i}} \Big] = \Big[ K^{\alpha,\beta}(z_i,z_j) \Big] .
\]
This is easily seen to be a continuous function from $PS^2$ into $\mr{GL}(n)$
such that $\Phi_F$ is unitarily equivalent to the map
\[ \pi(f)(\alpha,\beta) = Q_{\alpha,\beta}^{1/2} D_f Q_{\alpha,\beta}^{-1/2} .\]

\begin{thm}
Let $n \ge 3$, and let $\{ z_2,\dots,z_n \}$ be a set of 
distinct non-zero points in $\bb D.$ 
There exists $r > 0$ so that if $|z_1| \le r$ and 
$F= \{z_1, z_2,\dots,z_n \}$, then 
$\tilde\Phi_F: H^\infty_1/\cl I^F \to \rm C(PS^2, \ff M_n)$ 
is not a complete isometry.
\end{thm}

\begin{proof} 
Assume otherwise. 
Then there would exist a sequence of points 
$\{z_1(m) \}_{m=1}^{\infty}$ tending to 0 such that for each 
set $F_m = \{z_1(m), z_2,\dots,z_n \}$, the homomorphism 
$\tilde\Phi_{F_m}: H^\infty_1/\cl I^{F_m} \to \rm C(PS^2, \ff M_n)$ 
is a complete isometry. 
Set $F = \{0,z_2,\dots,z_n\}$.

Let $Q_{\alpha, \beta}(m)$ denote the matrix functions defined above
for the set $F_m$; and let $Q_{\alpha, \beta}$ denote the
matrix function corresponding to $F$ as before.
By the continuity of the function $K^{\alpha, \beta}$, 
as $m$ tends to $\infty,$ the functions $Q_{\alpha, \beta}(m)$ 
converge uniformly to $Q_{\alpha, \beta}$.
Using continuity and compactness, one can see that there exists 
a $\delta >0$ such that $Q_{\alpha, \beta}(m) \ge \delta I_n,$ 
for all $(\alpha, \beta) \in PS^2$ and all $m\ge1$.

Suppose that $W_1,\dots,W_n$ are $k \times k$ matrices satisfying 
\[
 \Delta(\alpha,\beta) :=
 \Big[ (I_k - W_iW_j^*) K^{\alpha, \beta}(z_i,z_j) \Big] \ge 0
 \qforal (\alpha,\beta) \in PS^2 .
\]
We claim that there exist a sequence $\epsilon_m \to 0$ such that 
\[
 \Big[ ((1+ \epsilon_m) I_k - W_iW_j^*) K^{\alpha, \beta}(z_i(m),z_j(m)) \Big] \ge 0
\]
for all $(\alpha,\beta) \in PS^2$ and all $m \ge 1$.
To see this, note that the difference between $ \Delta(\alpha,\beta)$ and
$\Delta_m(\alpha,\beta) := \big[ (I_k - W_iW_j^*) K^{\alpha, \beta}(z_i(m),z_j(m)) \big]$
is a sequence of Hermitian valued functions converging uniformly to $0$.
Thus, we may chose positive scalars $\epsilon_m \to 0$ such that 
\[
 \Delta_m(\alpha,\beta) - \Delta(\alpha,\beta)
  \ge -\epsilon_m \delta I_k
  \ge - \epsilon_m \big[ K^{\alpha, \beta}(z_i(m),z_j(m)) \big] .
\] 
The claim follows.
 
This inequality implies that there exists functions 
$f_m \in \ff M_k(H^\infty_1)$ with $\|f_m\|_\infty \le 1 + \epsilon_m$
such that $f_m(z_i(m)) = W_i$ for $i=1,\dots,n$.
Taking a weak*-limit point of these functions yields a function 
$f \in \ff M_k(H^\infty_1)$ with $\|f_m\|_\infty \le 1$ 
satisfying $f(z_i) = W_i$ for $i=1,\dots,n$.

This proves that $\pi: H^\infty_1/\cl I^F \to \rm C(PS^2, \ff M_n)$ 
is a complete isometry, contrary to Corollary~\ref{C:not c.isom}.
\end{proof}

Theorem~\ref{T:fail matrix} is now an immediate consequence.

\section{A Modified Pseudo-Hyperbolic Metric} \label{S:metric}

 M. B. Abrahamse \cite{Ab} proves that, in a certain sense,
all of the kernels given by his parameters are necessary for his
Nevanlinna-Pick type result. Further results about the necessity of all of Abrahamse's kernel conditions can be found in \cite{FV}, \cite{Mc} and \cite{McP}.
We consider a similar problem for our kernels in this section and examine in detail the pseudo-hyperbolic metric on the disk induced by the algebra $H^{\infty}_1.$

A uniform algebra $\cl A$ on a topological space $X$ induces a 
metric on the space by setting 
\[ d_{\cl A}(x,y) = \sup \{ |f(x)|: f \in \cl A, \|f\| \le 1, f(y)=0 \} .\]
Note that $0 \le d_{\cl A}(x,y) \le 1$; and that $d_{\cl A}(x,y) >0$
when $x \ne y$ because $\cl A$ separates points.
It is an elementary exercise with M\"{o}bius maps to show that 
$d_{\cl A}(x,y) = d_{\cl A}(y,x)$, and that 
this distance is comparable to the usual metric induced on $X$
considered as a subset of $\cl A^*$.  
This second distance is given by $\|\delta_x - \delta_y\|$,
where $\delta_x$ denotes the point evaluation $\delta_x(f) =f(x)$.
Moreover, it is a standard exercise to show that
\[ d_{\cl A}(x,y) = \dfrac{2 \|\delta_x - \delta_y\|}{1 + \|\delta_x - \delta_y\|^2} .\]
  From this, one can deduce that
\[
 d_{\cl A}(x,z) \le 
 \dfrac{ d_{\cl A}(x,y) + d_{\cl A}(y,z) }{ 1 + d_{\cl A}(x,y) d_{\cl A}(y,z)} .
\]
The triangle inequality is evident now, as is the fact that $d_{\cl A}(x,y) < 1$
is an equivalence relation.  The equivalence classes are called Gleason parts.
See \cite{Bear} for this material and its consequences.

The best known example of this construction is the 
{\em pseudohyperbolic metric,} which is the metric on $\bb D$ 
induced by $H^{\infty}.$ 
This metric is given by the formula
$d_H(z,w) = \left|\dfrac{z-w}{1- \overline{w}z}\right|.$ 
Note that $d_H(z,0) = |z|.$

In this section, we compute the metric $d_1$ on the disk induced by
the uniform algebra $H^{\infty}_1.$ 
Consideration of this metric will show that, at least many of the kernels 
$K^{\alpha,\beta}$ are necessary even for interpolation on two points.
While we have no need here of what $d_1$ is on the rest of the
maximal ideal space of $H^\infty_1$, which coincides with that of $H^\infty$,
it is easy to see that the two metrics are equal except when both points lie
in the open disk $\bb D$.  Indeed, if $x$ lies in the corona,
then replacing $f$ by $z^2f$ has no impact on the supremum 
in the definition of $d_{\cl A}(x,y)$.

First we compute 
\begin{align*}
d_1(z,0) &= \sup \{ |f(z)|: f \in H^{\infty}_1, \|f\|_{\infty} \le 1, f(0)=0 \}\\
&= \sup \{ |z^2g(z)|: g \in H^{\infty}, \|g\|_{\infty} \le 1 \} = |z|^2 .
\end{align*}

More generally, we can compute $d_1(z,w)$ by using Theorem~1.1. 
We know that there exists $f \in H^{\infty}_1$ with $f(w)=0$,
$\|f\|_{\infty} \le 1$ and $f(z)= \lambda$ if and only if
\[
 \begin{bmatrix}
  K^{\alpha,\beta}(w,w) & K_{\alpha, \beta}(w,z) \\
  K^{\alpha,\beta}(z,w) & (1- |\lambda|^2)K_{\alpha, \beta}(z,z) 
 \end{bmatrix} 
 \ge 0 \qforal |\alpha|^2 + |\beta^2 = 1 .
\]
Since the diagonal entries are positive, this latter condition holds 
if and only if the determinant is non-negative; i.e.
\[
 |\lambda|^2 \le 1 - \frac{|K^{\alpha,\beta}(w,z)|^2}
 {K^{\alpha,\beta}(w,w) K^{\alpha,\beta}(z,z)} 
 \qforal |\alpha|^2 + |\beta^2 = 1 .
\]
Therefore
\[
 d_1(z,w)^2 = \min \Big\{ 1 - \frac{|K^{\alpha,\beta}(w,z)|^2}
 {K^{\alpha,\beta}(w,w) K^{\alpha,\beta}(z,z)}
 : |\alpha|^2+|\beta|^2 =1 \Big\}.
\]

When $w=0$ and  $z=re^{i \theta}$, this simplifies to 
\[
 d_1(z,0)^2 = \min \Big\{  \frac{r^4}{r^4 + (1-r^2) |\alpha + \beta z|^2} .
 : |\alpha|^2+|\beta|^2 =1 \Big\}.
\]
It is clear that the minimum occurs precisely when
\[
 \alpha =  \frac{e^{i(s+\theta)}}{\sqrt{1+r^2}} \qand  \beta= \frac{re^{is}}{\sqrt{1+r^2}} 
\] 
for any $e^{is} \in \bb T$; and equals $r^4$, in agreement with the
earlier calculation.

Note that as $z$ varies over the disk, the pair $(\alpha,\beta)$ where
the minimum is attained, up to multiplication by a complex scalar of
modulus one, exhausts all of the kernel functions $K_{\alpha, \beta}$ 
with $|\alpha| \ge 1/\sqrt 2$. 
Therefore, even for interpolation of two points, $\{0,z\}$,
all of the kernels for $|\alpha| \ge 1/\sqrt 2$,  
or at least a dense subset, are necessary to attain the appropriate 
value for the minimum as $z$ varies over the disk.

\section{$C^*$-envelopes for Two-point Interpolation} \label{S:two points}

As we have seen in Section~\ref{S:Cenv} and also in the work of 
\cite{McP} and \cite{So}, $C^*_e(\cl A/\cl I^F)$ can be quite difficult 
to understand and be quite a complicated C*-algebra when $\cl I^F$ 
is the ideal of functions vanishing at 3 or more points. 
In contrast, given any uniform algebra $\cl A$ on a compact 
Hausdorff space $X$ and any two point set $F$, 
$C^*_e(\cl A/\cl I^F) = \ff M_2$ or $C^*_e(\cl A/\cl I^F) = \bb C \oplus \bb C.$ 
This fact is referred to in \cite{McP}, but we shall make it a bit more 
explicit here by exhibiting the completely isometric representation 
of $\cl A/\cl I^F$ into $\ff M_2$.  This will highlight the relationship of this 
representation with the pseudo-hyperbolic metric introduced in the 
previous section.

Let us assume that $X$ is a compact, Hausdorff space and
$\cl A \subseteq \rm C(X)$ is a uniform algebra.
Fix a two point subset $F= \{x_1, x_2 \}$ of $X$, 
and let $\cl I^F$ denote the ideal of functions vanishing on $F$.

Pick any two functions, $f_1 ,f_2 \in \cl A,$ such that $f_i(x_j) = \delta_{i,j}$.
It is easy to see that in the quotient, $E_i =f_i + \cl I^F$ for $i=1,2$  satisfy
\[ E_1^2=E_1,\ \ E_2^2=E_2,\ \ E_1E_2=E_2E_1 = 0\ \AND \ E_1+E_2 = E ,\] 
where $E= 1+\cl I^F$ denotes the identity element of the quotient algebra.
That is, $\cl A/\cl I^F$ is a two-idempotent operator algebra in the sense of \cite{Pa}.
Moreover, as is discussed in \cite{Pa}, given $w_1,w_2 \in \bb C,$ 
\[ \|w_1E_1 + w_2 E_2 \| = \inf \{ \|f\|: f \in \cl A, f(x_i) =w_i, i=1,2 \} .\]
More generally, given $W_1,W_2 \in \ff M_k,$ we have that
\[
 \| W_1 \otimes E_1 + W_2 \otimes E_2 \| 
 = \inf \{ \|F\|: F \in \ff M_k(\cl A), F(x_i) =W_i, i=1,2 \}.
\]
In particular, we have that 
\[
  d_{\cl A}(x_1,x_2) = \sup \{ |w_1|: \|w_1E_1\| \le 1 \} = \|E_1\|^{-1} .
\]
By a similar argument, $d_{\cl A}(x_1,x_2) = \|E_2\|^{-1}$.
 
Let $\pi:\cl A/\cl I^F \to B(\cl H)$ be a completely isometric representation
of $\cl A/\cl I^F$ as operators on some Hilbert space.
Then $\pi(E_1)$ and $\pi(E_2)$ will be idempotent operators which sum to the identity.
Hence we may decompose $\cl H= \cl H_1 \oplus \cl H_2$,
where $\cl H_1$ is the range of operator $\pi(E_1)$. 
Writing $\pi(E_1)$ and $\pi(E_2)$ as operator matrices with respect to 
this decomposition, we see that there exists a bounded operator 
$B: \cl H_2 \to \cl H_1$ such that
\[
 \pi(E_1) = \begin{bmatrix} I_{\cl H_1} & B\\0 & 0 \end{bmatrix} \qand
 \pi(E_2) = \begin{bmatrix} 0 & -B\\0 & I_{\cl H_2} \end{bmatrix}.
\]
Since $d_\cl A(x_1,x_2)^{-2} = \|E_1\|^2 = (1 + \|B\|^2),$ we see that the norm of $B$ is determined by the pseudo-hyperbolic metric.

As $\pi$ is completely isometric representation, we have that
\begin{align*}
 \|W_1 \otimes E_1 + W_2 \otimes E_2 \| &= 
 \bigg\| \begin{bmatrix} 
  W_1 \otimes I_{\cl H_1} & (W_1 -W_2) \otimes B\\
  0 & W_2 \otimes I_{\cl H_2} 
 \end{bmatrix} \bigg\| \\[.5ex] &=
 \bigg\| \begin{bmatrix}
  W_1 & (W_1 - W_2)\|B\| \\ 0 & W_2
 \end{bmatrix} \bigg\| .
\end{align*}
The last equality follows by computing the norm of the middle term.

The above observations lead readily to the following.

\begin{thm} 
Let $X$ be a compact, Hausdorff space; 
and let $\cl A \subseteq \rm C(X)$ be a uniform algebra.
Let $F=\{x_1,x_2\}$ be a two element subset of $X$;
and set $b= (d_\cl A(x_1,x_2)^{-2} -1)^{1/2}$.
Then the representation $\pi:\cl A/\cl I^F \to \ff M_2$ defined by 
\[
  \pi(f+\cl I^F) = 
  \begin{bmatrix} f(x_1) & b \big( f(x_1) -f(x_2) \big)\\ 0 & f(x_2) \end{bmatrix}
\]
is completely isometric.
Consequently, 
\[
 C^*_e(\cl A/\cl I^F) = 
 \begin{cases}
  \ff M_2, & \text{when } d_\cl A(x_1,x_2) < 1\\ 
  \bb C \oplus \bb C & \text{when } d_\cl A(x_1,x_2) =1 
 \end{cases} .
\]
\end{thm}

In the case that $\cl A= H^\infty_1$ and $F = \{x_1, x_2\} \subset \bb D$, 
it is readily seen that $d_\cl A(x_1,x_2) < 1$. 
Therefore $C^*_e(H^\infty_1/\cl I^F) = \ff M_2.$

\section{Concluding Remarks and Open Problems} \label{S:final}

We have seen that for some sets of points, the matrix-valued version of the 
interpolation result (Theorem~\ref{thm1}) fails, in particular, 
if one of the points is sufficiently close to $0$. 

\begin{prob} \label{P1}
Given a finite subset  $F=\{z_1,\dots,z_n\}$ of $\bb D$ with $n\ge3$,
is the homomorphism $\tilde\Phi_F: H^\infty_1/\cl I^F \to \rm C(PS^2, \ff M_n)$ 
ever completely isometric? Is it, in fact, ever 2-isometeric?
\end{prob}

\begin{prob} For three distinct points of the form $\{ 0,z_2,z_3 \},$
find three explicit matrices for which interpolation fails?
Does the homomorphism $\tilde\Phi_F$ fail to be even two isometric?
\end{prob}

\begin{prob} Disprove for $f=(f_{i,j}) \in \ff M_k(L^\infty),$ the distance formula
\[ 
 \dist(f, \ff M_k(H^\infty_1)) =
 \sup_{|\alpha|^2+|\beta|^2=1}
 \|((I-P_{\alpha,\beta})M_{f_{i,j}} P_{\alpha,\beta})\|.
\]
by exhibiting a concrete function.
\end{prob}

Although we know what the C*-envelope of the quotient algebra 
is in the case where one of the points is $0$, we have not been able 
to determine the C*-envelope in any other cases. 
The natural follow-up to Problem~\ref{P1} is:

\begin{prob} For $n \ge 3$ distinct non-zero points in $\bb D,$ 
find the C*-envelope, $C^*_e(H^\infty_1/\cl I^F).$ 
Are its irreducible representations all of dimension $n+1$?
\end{prob}

In the discussion of the hyperbolic metric, we showed that a large
set of kernels are necessary to determine interpolation.  But we were
not able to show that all are required.

\begin{prob} 
Is a dense set of kernels, $K^{\alpha,\beta}$, (up
to multiplication by a complex number of modulus one), necessary to
attain the metric $d_1$? If not, is a dense set necessary for the
general interpolation problem?
\end{prob}

A famous result for the annulus due to Federov-Vinnikov \cite{FV} 
says that once one fixes the points $z_1,\dots,z_n$, 
then to determine interpolation, one does not need to consider the 
whole family of kernel functions parameterized by the circle, 
but in fact, for scalar interpolation, there are two points on the circle 
such that just these two kernel functions will give necessary and 
sufficient conditions for interpolation. See \cite{Mc} for another proof of this fact. 
Thus, although as the points vary over the annulus all kernels are 
needed, once the points are specified only two kernels are needed 
for scalar-valued interpolation. 
Later, McCullough\cite{Mc} proved that even when the points are 
specified, to obtain necessary and sufficient conditions for 
matrix-valued interpolation all kernel functions (that is a dense 
subset) are needed.
McCullough's result was refined somewhat by McCullough and 
the second author in \cite{McP}, where it was shown that for any 
finite set of three or more points $z_1,\dots,z_n$ in the annulus the 
C*-envelope of the quotient algebra was isomorphic to 
$\ff M_n \otimes \rm C(\bb T)$; and that the irreducible representations 
of this algebra were parameterized by the kernel functions. 
McCullough's result then followed from this computation of the C*-envelope.

These considerations motivate the following problem.

\begin{prob} 
Given $z_1,\dots,z_n$ distinct points in $\bb D$,
does there exist a finite subset $F$ of the complex unit two-sphere
such that, given any complex numbers, $w_1,\dots,w_n$, 
there exists $f \in H^{\infty}_1$ with $\|f\|_{\infty} \le A$  
and $f(z_i)=w_i$ for $i=1,\dots,n$ if and only if 
\[ \Big[ (A^2 - w_i \bar{w_j})K^{\alpha,\beta}(z_i,z_j) \Big] \]
is positive semidefinite for all $(\alpha,\beta) \in F$?
\end{prob}



\begin{thebibliography}{99}

\bibitem{Ab} M. B. Abrahamse, 
{\em The Pick Interpolation Theorem for Finitely Connected Domains,} 
Michigan Math.\ J. {\bf 26} (1979), 195--203.
    
\bibitem{AM3} J. Agler and J.E. McCarthy,     
\textit{Complete Nevanlinna-Pick kernels},
J. Funct.\ Anal.\ \textbf{175} (2000), 111--124.

\bibitem{AM} J. Agler and J.E. McCarthy, 
{\em Pick Interpolation and Hilbert Function Spaces,} 
Graduate Studies in Mathematics, \textbf{44}, 
American Mathematical Society, Providence, RI, 2002. 

\bibitem{AMc2} J. Agler and J.E. McCarthy,     
\textit{Distinguished varieties},
Acta Math.\ \textbf{194} (2005), 133--153.
    
\bibitem{Bear} H.S. Bear,
\textit{Lectures on Gleason parts },
Lecture Notes in Math.\ \textbf{121},
Springer-Verlag, New York, 1970.

\bibitem{FV} S.I. Fedorov and V.L. Vinnikov, 
{\em On the Nevanlinna-Pick interpolation in multiply connected domains,} (Russian) 
Dept.\ of Mathematics University of Auckland Report Series {\bf 325}, October 1995:
Zap.\ Nauchn.\ Sem.\ St.\ Petersburg  Otdel.\ Mat.\ Inst.\ Steklov (POMI) {\bf 254} (1998); 
Anal.\ Teor.\ Chisel.\ i Teor.\ Funkts.\ {\bf 15}, 5--27. 

\bibitem{He} H. Helson, 
{\em Harmonic Analysis,} 
Addison-Wesley  Publishing Company, 
Reading, MA,  1983.
  
\bibitem{Mc} S. McCullough, 
{\em Isometric representations of some quotients of $H^{\infty}$ of an annulus,} 
Integral Equations Operator Theory {\bf 39} (2001), 335--362.

\bibitem{McP} S. McCullough and V.I. Paulsen, 
{\em C*-envelopes and Interpolation Theory,} 
Indiana Univ.\ Math.\ J. {\bf 51} (2002), 479--505.
    
\bibitem{Pa} V.I. Paulsen, 
{\em Operator Algebras of Idempotents,}
J. Funct.\ Anal.\ {\bf 181} (2001), 209--226.    

\bibitem{Pau} V. Paulsen,
\textit{Completely bounded maps and operator algebras},
Cambridge Studies in Advanced Mathematics \textbf{78},
Cambridge University Press, Cambridge, 2002.

\bibitem{PS} V.I. Paulsen and D. Singh, 
{\em Modules Over Subalgebras of the Disk Algebra}, 
Indiana Univ.\ Math.\ J. \textbf{55} (2006), 1751--1766.

\bibitem{Sa} D. Sarason,
\textit{Generalized interpolation in $H^\infty$},
Trans.\ Amer.\ Math.\ Soc. \textbf{127} (1967), 179--203.

\bibitem{So} J. Solazzo, 
\textit{Interpolation and Computability}, 
PhD Thesis,  University of Houston, 2000.

\end{thebibliography}
\end{document}